\documentclass[11pt]{amsart}

\usepackage{amssymb,amsmath}
\usepackage{bbm}
\usepackage{a4wide}

\newtheorem{theorem}{Theorem}
\newtheorem{prop}{Proposition}
\newtheorem{lemma}{Lemma}
\newtheorem{coro}{Corollary}
\newtheorem{fact}{Fact}
 
{
\theoremstyle{definition}

\newtheorem{remark}{Remark}
\newtheorem{example}{Example}
}

\newcommand{\ts}{\hspace{0.5pt}}
\newcommand{\CC}{\mathbb{C}\ts}
\newcommand{\RR}{\mathbb{R}\ts}
\newcommand{\QQ}{\mathbb{Q}\ts}
\newcommand{\ZZ}{\mathbb{Z}}
\newcommand{\TT}{\mathbb{T}}

\newcommand{\FF}{\mathbb{F}}
\newcommand{\KK}{\mathbb{K}}

\newcommand{\one}{\mathbbm{1}}
\newcommand{\sph}{\mathbb{S}}

\newcommand{\cH}{\mathcal{H}}
\newcommand{\cF}{\mathcal{F}}
\newcommand{\cS}{\mathcal{S}}
\newcommand{\cR}{\mathcal{R}}

\DeclareMathOperator{\ord}{ord}

\begin{document}

\title[Reversing symmetry groups]
{The structure of reversing symmetry groups}

% \vspace*{5mm}
\author{Michael Baake}
\address{Fakult\"at f\"ur Mathematik, Universit\"at Bielefeld, 
Box 100131, 33501 Bielefeld, Germany}
\email{mbaake@math.uni-bielefeld.de}
\urladdr{http://www.math.uni-bielefeld.de/baake}

\author{John A.~G.~Roberts}
\address{School of Mathematics, 
University of New South Wales,
Sydney, NSW 2052, Australia}
\email{jag.roberts@unsw.edu.au}
\urladdr{http://www.maths.unsw.edu.au/\~{}jagr}

\begin{abstract} 
  We present some of the group theoretic properties of reversing
  symmetry groups, and classify their structure in simple cases that
  occur frequently in several well-known groups of dynamical systems.
\end{abstract}

\maketitle

\vspace*{-5mm}
\section{Introduction}

Let $X$ be some space, with automorphism group $G:=\text{Aut}(X)$.  An
element $L\in G$ is said to have a {\em symmetry}\/ if there exists an
automorphism $S\in G$ that satisfies $L\circ S = S\circ L$ or,
equivalently, $S \circ L \circ S^{-1} = L$, and a {\em reversing
symmetry}, or {\em reversor}, if there exists an automorphism $R\in G$
so that $ R \circ L \circ R^{-1} = L^{-1}$.  The set of symmetries is
non-empty (it certainly contains all powers of $L$) and forms a group,
the {\em symmetry group}\/ $\cS(L)$. On the other hand, the existence
{\em a priori}\/ of any reversing symmetries for a particular $L$ is
unclear. When $L$ has a reversing symmetry, $L$ is called {\em
reversible}, and {\em irreversible}\/ otherwise.  The set $\cR(L)$ of
all symmetries and reversing symmetries of $L$ is a group, too, called
the {\em reversing symmetry group}\/ \cite{lamb} of $L$ (see also
\cite{Goodson}).

The simultaneous consideration of ordinary and reversing symmetries of
reversible automorphisms (which may arise as the time-one maps of
reversible flows) is known to provide some powerful algebraic
insights. As the results of \cite{lamb,Goodson} illustrate, the
knowledge of $\cS (L)$ has several implications on the nature of
possible reversing symmetries in $\cR (L)$. The power of this group
theoretic setting has recently been realized in cases where one has
access to the structure of the symmetry group $\cS(L)$, as in the case
of toral automorphisms \cite{BRcat,BRtorus} (via Dirichlet's unit
theorem \cite[Ch.\ 15.5]{Hasse}) or polynomial automorphisms of the
plane \cite{RBstandard,BRpoly,GM,GM2} (via the classification of
Abelian subgroups according to \cite{Wright}).
 
In many cases of reversible automorphisms (and also in the analogous
continuous-time case of reversible flows), it is in fact found that
all its reversing symmetries $R$ are involutions or elements of small
even order. Whenever an involutory reversor exists, the automorphism
$L$ can be written as the composition of two involutions, e.g., $L
\circ R$ and $R$, or $R$ and $R \circ L$, an observation that goes
back to Birkhoff \cite{birk}.  References \cite{roqu} and \cite{LR}
include reviews of the properties and applications of reversible
automorphisms and flows.

The goal of the present paper is to analyze the general structure of
$\cR(L)$, distilling and extending some theoretical insights from
specific cases already considered in
\cite{BRcat,BRtorus,BRpoly,RBstandard}.  Most of the algebraic methods
we use below are standard. Nevertheless, in view of the applications
to dynamical systems, we try to make the text self-contained
as far as algebraic methods are concerned (also giving references for
further background material).

\section{Mathematical setting}

From now on, we shall work within a given group $G$, e.g., the
automorphism group of some space $X$. Elements of $G$ will be denoted
by $f,g$ etc., with $1$ being the neutral element.  Motivated by the
dynamical systems context, we define the following subgroups of
$G$. The {\em symmetry group}\/ of an element $f\in G$ is the
centralizer of this element within $G$, i.e.,
\begin{equation} \label{def-sym}
   \cS(f) \; := \; \text{cent}^{}_{G} (f) \; = \;
        \{ g\in G\mid f g = g f \} \; = \;
        \{ g\in G\mid g f g^{-1} = f \} \, .
\end{equation}
The {\em reversing symmetry group}\/ $\cR(f)$ is defined as
\begin{equation} \label{def-rev}
    \cR(f) \; := \; \{ h\in G\mid hfh^{-1} = f^{\pm 1} \} \, .
\end{equation}
% \end{defin}

There are well-known facts about the groups $\cS(f)$ and $\cR(f)$, not
all of which are easy to locate in the literature. In this section, we
recall and extend some results that are relevant to our later
discussion, providing short proofs.

\smallskip
Clearly, $\cR(f)$ is a subgroup of $G$ that contains $\cS(f)$, and
one has, compare \cite{lamb,BRcat}:
\begin{fact} \label{rev-as-extension}
   $\cS(f)$ is a normal subgroup of\/ $\cR(f)$, with the factor
   group\/ $\cR(f)/\cS(f)$ either being the trivial group or\/
   $C_2$, the cyclic group of order\/ $2$.
\end{fact}
\begin{proof}
If $\cR(f)=\cS(f)$, which happens if $f^2=1$ or if $f$ is
irreversible, the statement is trivial. So, assume that there is an
$r\in\cR(f)$ with $rfr^{-1}=f^{-1}\neq f$. Define a binary {\em
grading}\/ $\Sigma\!:\, \cR(f)\longrightarrow C_2 = (\{\pm 1\},\cdot)$
by $\Sigma(h):=\varepsilon$ when $hfh^{-1}=f^{\varepsilon}$.  This
grading is a group homomorphism with ${\rm ker}(\Sigma)=\cS(f)$
(whence $\cS(f)$ is a normal subgroup) and ${\rm im}(\Sigma)=C_2$
(whence $\cR(f)/\cS(f) \simeq C_2$), which establishes the claim.
The grading highlights the fact that the composition of two
reversors is a symmetry.
\end{proof}

The case that $f^2=1$ is not of particular interest, as this
always gives $\cR(f)=\cS(f)$, due to $f^{-1}=f$. So, from now 
on, we shall always assume
\begin{itemize}
\item The element $f\in G$ satisfies $f^2\neq 1$, i.e., it is
   neither $1$ nor an involution.
\end{itemize}
This has an immediate consequence \cite[Prop.~5]{lamb} that we shall
need a number of times below:
\begin{fact} \label{not-odd}
  If\/ $f$ is reversible, with\/ $f^2 \neq 1$, no
  reversor of\/ $f$ can be of odd order.
\end{fact}
\begin{proof}
If $r$ is any reversor of $f$, we have $rfr^{-1}=f^{-1}$. This implies
$rf^{-1}r^{-1}=f$ (using $1=rr^{-1}=rfr^{-1}\, rf^{-1}r^{-1} = f^{-1}
\, rf^{-1}r^{-1}$), hence also $r^{\ell}f r^{-\ell}=f^{(-1)^{\ell}}$.
Then, $r^{2m+1}=1$ would give $f=f^{-1}$, contradicting the
assumption.
\end{proof}

\begin{remark}
For most applications in dynamical systems, one is mainly interested
in the situation that $f$ is not of finite order, so that $\langle f
\rangle := \{ f^n\mid n\in\ZZ\} \simeq C_{\infty}$. In this case, the
reversing symmetry group $\cR(f)$ of \eqref{def-rev} can also be
formulated as
\begin{equation} \label{normal}
   \cR(f) \; = \; \mathrm{norm}^{}_{G} (\langle f \rangle) 
\end{equation}
because $h \langle f \rangle h^{-1} = \langle f \rangle$ is only
possible here if $f$ (as a generator of $\langle f \rangle$) is
conjugated into a generator, hence into either $f$ or $f^{-1}$.
Eq.~\eqref{normal} should then be compared with $\cS(f) =
\text{cent}^{}_{G} (\langle f \rangle)$, the latter being an obvious
reformulation of Eq.~\eqref{def-sym}.

If $f$ is of {\em finite}\/ order, $\mathrm{norm}^{}_{G} (\langle f
\rangle)$ contains $\cR (f)$ as a subgroup, but possibly further
elements, e.g., elements $h$ with $hfh^{-1}=f^2$. It might then be
advantageous, also in view of questions discussed in \cite{Goodson2},
to consider this extension.
\end{remark}

Let us recall another classic result on the order of reversing
symmetries, see \cite[Thm.~1.1.5]{Lamb-thesis}, with a considerably
simplified proof.
\begin{fact} \label{powerof2}
  Let\/ $f\in G$, with\/ $f^2\neq 1$, have a reversor\/ $r$ of finite
  order. Then, this order is\/ $\mathrm{ord} (r)=2^{\ell}
  (2m+1)$ for some\/ $\ell\ge 1$, and\/ $f$ also has a reversor\/ $r'$
  of order\/ $2^{\ell}$. The set of all reversors of\/ $f$,
  within\/ $G$, is thus given by\/ $r'\mathcal{S} (f)$.
\end{fact}
\begin{proof}
Clearly, the order of $r$ is even, by Fact~\ref{not-odd}, and hence of
the form stated. Define $r'=r^{2m+1}$, which is a reversor of $f$
because $2m+1$ is odd. Clearly, $r'$ has order $2^{\ell}$.
Fact~\ref{rev-as-extension} implies that we can use $1$ and $r'$ as
the coset representatives of $\cS (f)$ in $\cR (f)$, so that we get
$\mathcal{R} (f) = \mathcal{S} (f) \,\dot{\cup}\, r' \mathcal{S} (f)$.
\end{proof}

An important consequence of Fact~\ref{powerof2} is that we may restrict
the search for reversing symmetries to elements of order $2^\ell$,
$\ell\ge 1$, provided there is a finite order reversor at all.

As mentioned in the Introduction, a particularly frequent case in
applications is that of an involutory reversor. To formulate the
corresponding result \cite{BRcat}, we write $N\rtimes G$ for the
semi-direct product of the groups $N$ and $G$, with $N$ the normal
subgroup.
\begin{lemma} \label{semidirect}
   Let\/ $f\in G$ be a mapping with\/ $f^2\neq 1$ and symmetry group\/
   $\cS(f)$.  If\/ $f$ has an involutory reversor\/ $r$, the reversing
   symmetry group is\/ $\cR(f) = \cS(f)\rtimes C_2$, with\/
   $C_2=\langle r \rangle$.
\end{lemma}
\begin{proof}
Once again by Fact~\ref{rev-as-extension}, we know that $1$ and $r$
can be used as the coset representatives, i.e., $\cR(f) =
\cS(f)\,\dot{\cup}\, r\cS(f)$, all seen as subgroups or subsets of
$G$. As $r$ is an involution, $\cR (f)/\cS (f)\simeq
\langle r\rangle = C_2$, thus establishing the semi-direct product.
\end{proof}

Let us give an important example where {\em all}\/ reversors are
involutions, irrespective of the structure of $\mathcal{S}(f)$.

\begin{example} \cite{JRV,JR} \label{ex-one}
Let $E$ be an elliptic curve defined over a field $\KK$. It is
birationally conjugate to a Weierstra{\ss} form $W \!:\,
y^2=x^3+Ax+B$, with $A,B \in \KK$. It is well known, see
\cite{Silver}, that the points $W(\KK)$ on the curve $W$ (or $E(\KK)$
on the curve $E$) with coordinates in $\KK$ form an Abelian group with
associated group law ``$+$''. Let $G$ be the group of birational
transformations over $\KK$ that map $E$ to itself. Then, in the
typical case (i.e., when the curve does not permit complex
multiplication, which can be worked out using the so-called
$j$-invariant), $G$ has the form
\begin{equation} \label{translate}
     G \; \simeq \; \mathcal{T} \rtimes \{\pm\one\},
\end{equation}
where $\mathcal{T}$ is the group of translations on $W$, $P\mapsto
P+\varOmega$ with $\varOmega\in W(\KK)$, and $\pm\one$ stands for
$P\mapsto\pm P$. Clearly, any $f\in\mathcal{T}$, with $f^2\neq\one$,
has $\mathcal{S} (f) = \mathcal{T}\simeq W(\KK)$, i.e., $\cS (f)$ is
Abelian. The reversors of $f\in G$ are always involutions, $P\mapsto
-P +S$, for some $S\in W(\KK)$.

The structure of $W(\KK)$, and hence of $\mathcal{S} (f)$, is quite
general, depending on the field $\KK$. In particular, $W(\RR)$ is a
one-dimensional compact Lie group, while $W(\CC)\simeq \TT^2$, the
$2$-torus. Moreover, $W(\QQ)$ is a finitely generated Abelian group,
hence, by \cite[Thm.~I.8.5]{Lang},
\[
     W(\QQ) \; \simeq \; F\times (C_{\infty})^{r^{}_W} ,
\]
where $r^{}_W$ is the rank of the curve and $F$ is the finite
torsion group.
\end{example}

Looking more closely at this example, one realizes that the extra
structure, in comparison to Lemma~\ref{semidirect}, is that any
involutory reversor $r$ of $f$ actually conjugates {\em all}\/
elements of the group $\mathcal{S} (f)$ into their inverses, not just
$f$. This is a situation that is not a priori restricted to a
translation group structure as in \eqref{translate}.  An important
part of the semi-direct product structure $\mathcal{R} (f)=\mathcal{S}
(f)\rtimes \langle r\rangle$ in Lemma~$\ref{semidirect}$ is the
induced automorphism $\sigma$ on the normal subgroup,
\begin{equation}  \label{induced}
       \sigma(g) \; := \; rgr^{-1} \; = \; rgr, 
\end{equation}
for all $g\in\mathcal{S}(f)$. Given an involutory reversor $r$, the
interplay between the nature of $\sigma$ and the structure of
$\mathcal{S} (f)$ can be used effectively to determine the detailed
structure of the group $\mathcal{R} (f)$, as we shall see below in
Theorem~\ref{two-infty}. In particular, given an involutory reversor
$r$, the order of any other reversor $rg$ with $g \in \mathcal{S}
(f)$, necessarily of this form by Fact~\ref{rev-as-extension} and of
even order by Fact~\ref{not-odd}, follows from the equations
\begin{equation} \label{sigmaorder}
      (rg)^{2k} \; = \; (rgr^{-1}g)^k 
      \; = \;  (\sigma(g)\, g)^k \, , 
      \quad\mbox{for } k \in \ZZ.
\end{equation}
This has the following simple consequences (compare Example~\ref{ex-one}):
\begin{prop}\label{oneforall}
  Consider the element\/ $f\in G$ of Lemma~$\ref{semidirect}$ with an
  involutory reversor\/ $r$, an element\/ $g\in\cS (f)$ and\/
  $\sigma$ as in\/ $\eqref{induced}$. Then, one has:
\begin{itemize}
  \item[\rm (1)] $\sigma(g)=g^{-1}$ iff the reversor $rg$ is an involution.
  Consequently, $r$ is a simultaneous reversor for all elements of the
  group\/ $\mathcal{S} (f)$ iff all reversors of\/ $f$ are involutions;
  \smallskip
\item[\rm (2)] any finite order reversor of\/ $f$ must have order\/
  $2\ell$, where\/ $\ell$ is the order of some finite order symmetry
  of\/ $f$. So, if no non-trivial symmetry of finite order exists,
  there can only be reversors that are involutions or of infinite
  order. In this case, if\/ $\sigma(g')\ne (g')^{-1}$ for some
  $g'\in\mathcal{S}(f)$, the reversor\/ $rg'$ is of infinite order.
  \end{itemize}
\end{prop}
\begin{proof}
The first claim is obvious from Eq.~\eqref{sigmaorder}, used with
$k=1$. The second claim follows from the observation that
$(rg)^2=\sigma(g)\,g$ is a symmetry, hence either trivial or
not of finite order under the assumptions made.
\end{proof}

\begin{coro} \label{Abelian}
If all reversors of\/ $f$ are involutions, the symmetry group\/
$\cS (f)$ is Abelian.
\end{coro}
\begin{proof}
By part (1) of Proposition~\ref{oneforall}, a reversor $r$ of $f$
is a simultaneous reversor for all elements of $\cS (f)$. So, if
$a,b\in\cS (f)$, we have $rar=a^{-1}$, $rbr=b^{-1}$ and
$r (ab) r= (ab)^{-1}$. Consequently,
\[
   a^{-1} b^{-1} \; = \; r\ts ab\ts r \; = \; (ab)^{-1}
   \; = \; b^{-1} a^{-1}
\]
which gives $ab = ba$.
\end{proof}

That the converse of Corollary~\ref{Abelian} is not true is
illustrated below in Theorem~\ref{two-infty} and the associated examples.

\section{Implications from the symmetry group} \label{impli}

To further explore the group theoretic concequences, let us recall the
concept of a group extension, compare \cite[Thm.~15.3.1]{Hall} or
\cite[Sec.~I.14]{Huppert}. Fact~\ref{rev-as-extension} shows that we
need to look at cyclic $C_2$-extensions of the symmetry group, but not
all such extensions will give reversing symmetry groups. It is thus a
natural task to select and classify those that do.  We now present
first steps in this direction, building on previous work by various
authors \cite{lamb,Goodson,LR}.

The main point in using the group theoretic setting comes from the
consequences of the structure of $\cS(f)$ to that of $\cR(f)$. 
Classifying the structure of (non-trivial) reversing symmetry groups
then means:
\begin{enumerate}
\item Start from groups of the form $N={\rm cent}_G (f)$, 
      for some $f$ with $f^2\neq 1$;
\item Search for an $h\in G\setminus N$ with $hfh^{-1}=f^{-1}$;
\item Classify $H=N\,\dot{\cup}\,hN$, a $C_2$-extension of $N$, 
      according to its group structure.
\end{enumerate}
Our point here is that such a classification is a purely group
theoretic problem. In concrete examples and applications,
special conditions can then lead to further restrictions. 

\begin{remark}
It can also become meaningful, or even necessary, to consider the
equation $hfh^{-1}=f^{-1}$ only {\em up to symmetries}, i.e., to look
for solutions of $hfh^{-1}=s f^{-1}$ with $s\in\cS(f)$. If $s$ is of
finite order, some power of $f$, $f^k$ say, is reversible in the usual
sense. This is the basic mechanism of reversing $k$-symmetries,
compare \cite{Lamb-thesis,LQ} for details. A similar remark applies to
$k$-symmetries in comparison to ordinary symmetries. Conversely, a
reversible element $f\in G$ might have a {\em root}\/ in $G$ that is
not reversible itself, but satisfies such a more general equation. We
shall meet this situation below, in part (3) of Theorem
\ref{two-infty}.
\end{remark}

Let us continue with a general observation, which is a rather direct
consequence of a result of Goodson, see \cite[Prop.~2]{Goodson} and
the generalization mentioned afterwards, and \cite[Fact~11]{BRtorus}.
\begin{prop} \label{rev-order}
  Let\/ $f\in G$ be an element of infinite order, and assume that\/
  $\cS(f)= \cF \times \langle g \rangle$ where\/ $\cF$ is some finite
  group of order\/ $N\ge 1$ $($not necessarily Abelian$\ts)$, and\/
  $g$ is some generator $($then necessarily of infinite order$\ts)$.
  If\/ $r$ is a reversor of\/ $f$, $r$ is an element of\/ {\em finite}
  order. Its order is even and divides\/ $2N$.
\end{prop}
\begin{proof}
If $r$ is a reversor, $r^2$ is a symmetry, hence $r^2=s g^m$, for some
$s\in\cF$ and some integer $m$. Note that, due to the assumption of
the direct product structure, we always have $s g = g s$, even if
$\cF$ itself is not Abelian. Since the group $\cF$ is finite and of
order $N$, we know that $s^n=1$ for some $n\ge 1$ that divides
$N$. This implies that $r^{2n}=g^{mn}$.

As $f$ is not of finite order, but clearly an element of $\cS(f)$, we
may assume $f^N=g^k$ for some (positive) integer $k$ without loss of
generality, modifying the argument just used (in particular, $k\neq
0$, while $k>0$ might require to replace $g$ by $g^{-1}$).

Since $rf = f^{-1}r$ by assumption (hence also $rf^{\ell} =
f^{-\ell}r$, for all $\ell\in\ZZ$), we choose $\ell = mnN$ and obtain
$r g^{kmn} = g^{-kmn}r$. Since $g^{kmn}=r^{2nk}$, this implies $r\,
r^{2nk} = r^{-2nk} r$ and thus $r^{4nk}=1$, i.e., $r$ is of finite
order. Since $r^{2n}=g^{mn}$, this is only possible for $mn=0$, hence
$m=0$. This implies $r^{2n}=1$, so the order of $r$ divides $2N$. If
$f$ is not of finite order, it is not an involution, and $r$ can then
not be of odd order by Fact~\ref{not-odd} (hence also $r\neq 1$).
\end{proof}

\begin{remark}
An alternative way to state the result of Proposition~\ref{rev-order}
is the following. If $f$ is an element of infinite order, such that
the factor group $\mathcal{S}(f)/\langle f \rangle$ is finite, any
reversor $r$ of $f$ must be of finite order. In particular, $r^{2k}=1$
for some integer $k\ge 1$ that divides the order of the factor group.
\end{remark}

\begin{theorem} \label{one-infty}
   Let\/ $f\in G$ be an element of infinite order, with\/
   $\cS(f)\simeq C_{\infty}$. If\/ $f$ is reversible, one has\/
   $\cR(f) = \cS(f)\rtimes C_2 \simeq D_{\infty}$, and all reversors
   of\/ $f$ are involutions.
\end{theorem}
\begin{proof}
If $\cS(f)=\langle g \rangle\simeq C_{\infty}$, we must have $f=g^m$
for some $0\neq m\in\ZZ$. Let $r$ be any reversor of $f$, which must
then be an involution by Proposition~\ref{rev-order}. This gives the
general structure of $\mathcal{R}(f)$ as a semi-direct product by
Lemma~\ref{semidirect}.

In view of Proposition~\ref{oneforall}, we now have to look at
$\sigma (g)=rgr^{-1}$. By the previous argument, all reversors of
$f$ are involutions, hence necessarily $rgr^{-1}=g^{-1}$, and
$\mathcal{R} (f) = \langle g\rangle\rtimes\langle r\rangle
\simeq D_{\infty}$ is clear.
\end{proof}

The situation of Theorem~\ref{one-infty} looks rather special, but
actually occurs in some important dynamical contexts.

\begin{example}
Let $G$ be the space group of the integer lattice $\ZZ$ in dimension
one, which is $G=\ZZ\rtimes \mathrm{O}(1)=\ZZ\rtimes \{\pm 1\}$,
compare \cite{CM}. So, $G$ contains all Euclidean transformations that
map $\ZZ$ onto itself, and it has the structure \eqref{translate} with
$\mathcal{T}=\{T_m\!:\, x\mapsto x+m\mid m\in\ZZ\}\simeq\ZZ$.  Now,
take $f=T_n$ (with $n\neq 0$) as our mapping, the $n$-fold shift. This
is a standard mapping considered in symbolic dynamics, compare
\cite{LM}, and $G$ is a very natural group to embed it in. Clearly,
$\mathcal{S}(T_n)=\mathrm{cent}^{}_{G} (T_n) = \mathcal{T} \simeq
C_{\infty}$, while the map $x\mapsto -x$ is an involutory reversor for
$T_n$, noting that $(T_n)^{-1}=T_{-n}$.  Consequently, we have
$\mathcal{R}(T_n)\simeq C_{\infty}\rtimes C_2$, as in
Theorem~\ref{one-infty}.

Note that all involutions in $G$ are of the form $x\mapsto -x + m$,
with $m\in\ZZ$, and are always conjugate, within $G$, to either 
$x\mapsto -x$ or $x\mapsto -x+1$. The latter are conjugate via a
half-integer shift, hence not within $G$, but within some larger
group.
\end{example}

\begin{remark}
In the previous example, we could replace $\ZZ$ by $\QQ$ or $\RR$,
with obvious changes to the symmetry group, though the latter is no
longer isomorphic with $C_{\infty}$. Also, if $\varGamma$ is the
generic lattice in $\RR^d$, its space group \cite{CM} is $G=\varGamma
\rtimes\{\pm\one\}$, as inversion is then the only isometry of the
lattice. The previous example can now easily be extended to an
arbitrary translation $f\!:\, x\mapsto x+a$ with $0\neq
a\in\varGamma$.
\end{remark}

\begin{example} 
Consider $G=\mathrm{PGL}(2,\ZZ)$, the group of integer matrices with
determinant $\pm 1$, identified up to an overall sign, i.e.,
$\mathrm{PGL}(2,\ZZ) = \mathrm{GL}(2,\ZZ)/\{\pm \one\}$. By
Dirichlet's unit theorem, one can show \cite{BRtorus} that, if
$M\in\mathrm{PGL}(2,\ZZ)$ is not of finite order, its symmetry group
is $\mathcal{S}(M)=\mathrm{cent}_G (M)\simeq C_{\infty}$. A concrete 
example, even with $C_{\infty}=\langle M\rangle$, is
\begin{equation} \label{fibocat}
    M \, = \, \left[ \begin{matrix} 0 & 1 \\ 1 & 1 \end{matrix}
    \right]  , \quad \mbox{with reversors }\;
    R \, = \, \left[ \begin{matrix} 1 & 0 \\ 1 & -1 \end{matrix}
    \right] \mbox{ and }\;
    R' \, = \, R M \, = \,
    \left[ \begin{matrix} 0 & -1 \\ 1 & 0 \end{matrix}
    \right]  ,
\end{equation}
where we write $[M]$ for a matrix up to overall sign. Note that
both $R$ and $R'$ are involutions in $\mathrm{PGL}(2,\ZZ)$, but
they are not conjugate within this group.

In other examples, $M$ need not be a generator of $C_{\infty}$, as, in
general, such a matrix can have roots in $G$. Note that the spectrum
of $M$ in \eqref{fibocat} is only self-reciprocal up to an overall
factor of $-1$, whence $M$ is not reversible in $\mathrm{GL}(2,\ZZ)$,
though its square is (see below).
\end{example}

\begin{remark}
Two groups of dynamical systems isomorphic to $\mathrm{PGL}(2,\ZZ)$
are:
\begin{enumerate}
\item The group of $3$-dimensional invertible polynomial
  maps which preserve the Fricke-Vogt invariant
\[
    I(x,y,z) \; = \; x^2 + y^2 + z^2 - 2 xyz - 1
\]
  and fix the point $(1,1,1)$. Corresponding to $M$, $R$ and $R'$ 
above are, respectively, the Fibonacci trace map
$(x,y,z) \mapsto (y,z,2yz-x)$ and its reversors
$(x,y,z)\mapsto (z,y,x)$ and $(x,y,z)\mapsto
(2yz-x,z,y)$, see \cite{RBtrace} and references
therein for details.
\item The group of homeomorphisms of the $2$-sphere $\sph^2$ that are
induced by quotienting the action of a $\mathrm{GL}(2,\ZZ)$ matrix on
$\TT^2$ by the reflection in the origin. Whenever the
$\mathrm{GL}(2,\ZZ)$ matrix is hyperbolic, this yields a so-called
pseudo-Anosov map of $\sph^2$, see \cite{C,LMK} for details.
\end{enumerate}
\end{remark}

\begin{theorem}  \label{two-infty}
   If\/ $f\in G$ is a reversible element of infinite order  
   with symmetry group\/ $\cS(f)\simeq C_2 \times C_{\infty}$,
   all reversors must be involutions or elements of order\/ $4$.
   In particular, one finds precisely one of the following three
   situations.
\begin{itemize}
\item[{\rm (1)}] $\cR(f)\simeq C_2\times D_{\infty}$, iff all reversors
    of\/ $f$ are involutions. \smallskip
\item[{\rm (2)}] $\cR(f)\simeq C_{\infty}\rtimes C_4$, iff all reversors
    of\/ $f$ are elements of order\/ $4$. \smallskip
\item[{\rm (3)}] $\cR(f)\simeq \big(C_2\times C_{\infty}\big)
    \rtimes C_2$, iff there are reversors both of order\/ $2$ and\/ $4$.
\end{itemize}
\end{theorem}
\begin{proof}
By assumption and Proposition~\ref{rev-order}, we know that a reversor
in this setting must be an involution or an order $4$ element.

Let $\cS(f)\simeq C_2\times C_{\infty}$ with $C_{\infty}=\langle g
\rangle$ and an involutory symmetry $s$, which is then unique by the
structure of the group. If the reversor $r$ is an involution, one has
$\cR(f)\simeq \cS(f)\rtimes C_2$.  Since $rsr$ is also an involutory
symmetry, we get $rsr=s$ by uniqueness, and $r$ and $s$ commute. Since
$r\neq s$, this gives $\cR(f) \simeq (C_2\times C_{\infty})\rtimes
C_2$, with either $rgr^{-1}=g^{-1}$ (then simplifying to $\cR(f)\simeq
C_2\times D_\infty$) or $rgr^{-1} = s g^{-1}$ (in which case $f$ must
be an even power of $g$). Note that, in the latter case, $\varrho =
gr$ is an element of order $4$, and a reversor of $f$.

If $\cR(f)=C_2\times D_{\infty}$, we are in the situation of part (1)
of Proposition~\ref{oneforall}, as $r$ conjugates all generators of
$\mathcal{S} (f)$ into their inverses. Consequently, all reversors of
$f$ are involutions then.

If $f$ has a reversor $r$ of order $4$, $r^2=s$ is the unique
involutory symmetry of $f$, and $f=r^{2\varepsilon} g^m$ for
$\varepsilon\in\{0,1\}$ and some integer $m\ne 0$.  In particular,
$r^2$ and $g$ commute, and $rgr^{-1}$ is a symmetry of $f$, so that
$rgr^{-1} = r^{2k} g^{\ell}$ for $k\in\{0,1\}$ and some
$\ell\in\ZZ$. Clearly, in view of $rfr^{-1} = f^{-1}$, this forces
$\ell=-1$.

If $k=0$, $r$ is also a reversor of $g$, and we have $\cR(f)\simeq
C_{\infty}\rtimes C_4$. Reversors are of the form $r g^n$ or $r^3
g^n$, all of which have order $4$. This is the only case for $m$ odd,
while for $m$ even also $k=1$ is possible, i.e., $rgr^{-1} = r^{2}
g^{-1}$. This gives a group with the presentation
\[
   \cR(f) \; = \; \langle r,g \mid r^4 = 1 , \; 
                  r g^{\pm 1} = g^{\mp 1} r^{-1} \rangle
\] 
which is an index $2$ extension of $\cS(f)\simeq C_2\times C_{\infty}$, 
but does not look like a simple semi-direct product. However, 
$\rho = g^{-1} r$ is an involution that satisfies
$\rho g \rho = r^2 g^{-1}$, and it is a reversor for $f$.
This brings us back to $\cR(f)\simeq (C_2\times C_{\infty})\rtimes
C_2$, where the outer $C_2$ is generated by $\rho$.

This chain of arguments shows that the $3$ cases of the theorem are
both (algebraically) possible and exhaustive.
\end{proof}

\begin{remark}
Note that the meaning of the group $\big(C_2\times C_{\infty}\big)
\rtimes C_2 = \big( \langle s \rangle \times \langle g \rangle \big)
\rtimes \langle \rho \rangle$ in case (3) of Theorem~\ref{two-infty}
includes the induced automorphism $\rho g \rho^{-1} = s g^{-1}$. This
is the key difference to case (1), where a different induced
automorphism permits the simplification shown.
\end{remark}

Examples of all three cases of Theorem~\ref{two-infty} appear among
hyperbolic toral automorphisms (or cat maps) and polynomial
automorphisms of the plane:

\begin{example} \cite{BRcat,BRtorus,W}
Elements $M$ of the matrix group $\mathrm{GL}(2,\ZZ)$ that are not of
finite order, including the hyperbolic ones, have
$\mathcal{S}(M)\simeq C_2 \times C_{\infty}$, where $C_2 =
\{\pm\one\}$. Reversible elements $M$, and associated reversors $R$
(with subscripts indicating their order), which illustrate each case
of Theorem~\ref{two-infty} are:
\[ \begin{array}{cll}
   (1): & M = \begin{pmatrix} 1 & 2 \\ 1 & 3 \end{pmatrix}, 
        & R_2 = \begin{pmatrix} 1 & 0 \\ 1 & -1 \end{pmatrix}; \\[5mm]
   (2): & M = \begin{pmatrix} 5 & 7 \\ 7 & 10 \end{pmatrix}, 
        & R_4 = \begin{pmatrix} 0 & -1 \\ 1 & 0 \end{pmatrix}; \\[5mm]
   (3): & M = \begin{pmatrix} 1 & 1 \\ 1 & 2 \end{pmatrix}, 
        & R_2 = \begin{pmatrix} 1 & 0 \\ 1 & -1 \end{pmatrix}
          \mbox{ and }\;
          R_4 = \begin{pmatrix} 0 & -1 \\ 1 & 0 \end{pmatrix}.
\end{array}
\]
Note that the third case is closely related to the previous
$\mathrm{PGL}(2,\ZZ)$ matrix in Eq.~\eqref{fibocat}.
\end{example}

\begin{example} \cite{BRpoly,GM2,RBstandard}
Consider the case that $G$ is the group of planar polynomial
automorphisms with coefficients in the field $\KK$, i.e., polynomial
transformations $x'=P(x,y)$, $y'=Q(x,y)$ that have a polynomial
inverse (for ease of notation, we use $(x',y')$ for the image points),
see \cite{Essen} for general background material. Utilizing the
classical result that $G$ is an amalgamated free product of two
groups, consequently giving knowledge of the Abelian subgroups within
$G$ \cite{Wright}, it can be shown that $\mathcal{S} (f)$ is
isomorphic to either $C_{\infty}$ or $C_2\times C_{\infty}$, when
$\KK\in\{\QQ,\RR\}$ and $f$ is dynamically non-trivial. The latter
property means, in the language of \cite{BRpoly}, that $f$ is a
so-called CR element of $G$, hence neither conjugate to an affine nor
to an elementary mapping in $G$.

Reversible elements (for $G$ with $\KK\in\{\QQ,\RR\}$) illustrating
each case of Theorem~\ref{two-infty} are:
\begin{itemize}
\item[(1)]  $f\! : \; x' = x + p(y) \, , \; y' = y + q(x')\,$,
            with $p\neq q$ odd polynomials;  \\[1mm]
            $s \! : \; x' = -x \, , \; y' = -y$, an involution; \\[1mm]
            $r\! : \; x' = -x-p(y)\, , \; y' = y$, an involution; \\[1mm]
            $\mathcal{S} (f) = \mathrm{cent}_G (f) =
            \langle s\rangle\times\langle f\rangle
            \simeq C_2 \times C_{\infty}$; \\[1mm]
            $\mathcal{R} (f) = \langle s\rangle\times
            \big( \langle f\rangle\rtimes
            \langle r\rangle\big)
            \simeq C_2 \times D_{\infty}$. \smallskip
\item[(2)]  $f\! : \; x' = -x+y^3 \, , \; y' = -y-(x')^3\,$, \\[1mm]
            involutory symmetry $s$ as in case (1); \\[1mm]
            $r\! : \; x' = -y\, , \; y' = x$, an order $4$ reversor,
            with $r^2=s$;\\[1mm]
            $\mathcal{S} (f) = \mathrm{cent}_G (f) =
            \langle s\rangle\times\langle f\rangle
            \simeq C_2 \times C_{\infty}$; \\[1mm]
            $\mathcal{R} (f) = \langle f\rangle\rtimes
            \langle r\rangle
            \simeq C_{\infty}\rtimes C_4$. \smallskip
\item[(3)]  $f\! : \; x' = x + p(y) \, , \; y' = y + p(x')\,$,
            with $p$ an odd polynomial.  \\[1mm]
            involutory symmetry $s$ and reversor $r$ as in case (1); \\[1mm]
            $t \! : \; x' = y \, , \; y' = x+p(y)\,$,
            so that $f=t^2$; \\[1mm]
            $\mathcal{S} (f) = \mathrm{cent}_G (f) =
            \langle s\rangle\times\langle t\rangle
            \simeq C_{2} \times C_{\infty}$; \\[1mm]
            $\mathcal{R} (f) = \big( \langle s\rangle\times
            \langle t\rangle\big)\rtimes\langle r\rangle
            \simeq\big(C_2\times C_{\infty}\big)\rtimes C_2$.
\end{itemize}
In cases (1) and (2), $f$ has no root in $G$. In case (3), $r'=tr$ is
an order $4$ reversor of $f$.

Note that \cite{RV} provides a test for reversibility within
the group $G$, when a reduction of the polynomial maps to
finite fields is possible.

\end{example}

Following on from Theorem~\ref{two-infty}, it would be nice to have
some similarly simple classification of the group structure of
$\cR(f)$ for more general symmetry groups $\cS(f)$.  However, things
quickly become more involved, in particular if $f$ possesses roots in
$G$, which is a situation frequently met in practice.  If $f$ has no
roots in $G$, one can go further as follows.

\begin{theorem} \label{prime-case}
   Let\/ $\cS(f)\simeq C_p \times C_{\infty}$ with\/ $C_{\infty}=
   \langle f\rangle$ and\/ $C_p=\langle h \rangle$, $p$ an odd prime.
   If\/ $f$ is reversible, there are always involutory reversors, and
   one meets precisely one of the following two situations.
\begin{itemize}
\item[{\rm (1)}] $\cR(f)=C_{\infty} \rtimes C_{2p} = \langle f \rangle
   \rtimes \langle r\rangle$, with\/ $rfr^{-1}=f^{-1}$ and\/ $h=r^2$,
   iff a reversor $r$ of order $2p$ exists;
\item[{\rm (2)}] $\cR(f)=(C_{p}\times C_{\infty})\rtimes C_2$,
   with\/ $C_2=\langle r \rangle$ and\/ $rhr=h^{-1}$, iff
   all reversors are involutions.
\end{itemize}
\end{theorem}
\begin{proof}
By Proposition~\ref{rev-order}, any reversor $r$ of $f$ must be of
even order that divides $2p$, so either $\ord(r)=2$ or $\ord(r)=2p$
because $p$ is prime. In the latter case, in line with
Fact~\ref{powerof2}, $r^p$ is an involutory reversor, so that $\cR(f)
= (C_p \times C_{\infty})\rtimes C_2$ in both cases, by
Lemma~\ref{semidirect}.

We can thus focus on the equation $qfq^{-1}=f^{-1}$ with $q^2=1$, and
consider the possible automorphisms induced by $q$ on $\cS(f)$.
Clearly, $qhq^{-1}=qhq$ is a symmetry of $f$ of order $p$, so that
$qhq = h^m$ for some $1\le m\le p-1$. Since $h=q^2 h q^2 = q h^m q =
h^{m^2}$, we must have $m^2 = 1$ (mod $p$).  Since $p$ is a prime (and
$\FF_p$ thus a finite field), this congruence has precisely two
solutions. These are $m=\pm 1$ (mod $p$), either giving $qhq = h$ or
$qhq=h^{-1}$.

In the first case, $h$ and $q$ commute, and $(qh)^k = q^k h^k$.
This shows that $r=qh$, which is also a reversor, has order $2p$,
and $r^p=q$. Consequently, the reversing symmetry group becomes
$\cR(f)=C_{\infty}\rtimes C_{2p}$ with $C_{2p}=\langle r \rangle$
and $r^2=h$.

In the second case, $q$ is a reversor also for the finite order
element $h$, and the structure of $\cR(f)$ is as claimed. This brings
us back to the situation of part (1) of Proposition~\ref{oneforall},
whence all reversors of $f$ are involutions.
\end{proof}

\begin{remark}
If one considers $\cS(f)=C_n \times C_{\infty}=\langle h\rangle
\times\langle f \rangle$ for $n>1$ not a prime, things quickly become
more complicated.  In the case that $n$ is odd, any reversor $r$ must
have order $2\ell$ for some $\ell | n$, by
Proposition~\ref{rev-order}.  But then, $s=r^{\ell}$ is an involutory
reversor, and we can again restrict ourselves to looking at the
equation $sfs=f^{-1}$ and the induced automorphism on $\cS(f)$. In
this case, $shs=h^m$ for some $m\in\{1\le k\le n\mid \gcd(k,n)=1\}$,
subject to the additional requirement that $m^2 = 1$ (mod $n$). This
equation always has the solutions $m = \pm 1$ (mod $n$). They are the
only ones for $n=p^k$ with $k\ge 1$ and $p$ an odd prime, while more
solutions exist otherwise, e.g., $n=15$ permits $m=\pm 1$ and $m=\pm
4$. The number of solutions is $2^a$, with $a\ge 1$ the number of
distinct prime divisors of $n$, see \cite[Ch.\ 6.3]{Hasse}.  The
result of Theorem~\ref{prime-case} has to be extended accordingly.

If $n$ is even, such extra solutions may exist as well (e.g., $n=8$
permits $m=\pm 1$ and $m=\pm 3$, while $n=12$ is compatible with
$m=\pm 1$ and $m=\pm 5$). Here, if we write $n=2^{k+1} (2\ell+1)$ with
$k\ge 0$, the number of solutions is $2^{a+\min\{k,2\}}$, with $a\ge
0$ the number of distinct prime divisors of $2\ell +1$, compare
\cite[Ch.\ 6.3]{Hasse}. In general, it is no longer true that at least
one involutory reversor exists, as we already saw in case (2) of
Theorem~\ref{two-infty}.
\end{remark}

Next, let us take a closer look at a case where an additional symmetry
of infinite order exists. This is motivated both by the structure of
(projective) toral automorphisms in dimensions $d > 2$, see
\cite{BRtorus}, and by other examples from algebraic dynamics, see
\cite{Klaus} and references therein for an orientation.
\begin{theorem} \label{inf-inf}
   Let\/ $\cS(f)=\langle t\rangle \times\langle g\rangle\simeq
   C_{\infty}\times C_{\infty}$, with\/ $f=g^n$ for some integer\/
   $n\neq 0$. Then, any reversor\/ $r$ of\/ $f$ is either an
   involution $($hence giving\/ $\cR(f)$ as in
   Fact~$\ref{semidirect})$ or it is not of finite order.
 
   Moreover, $r$ is also a reversor for\/ $g$, and either\/
   $\sigma (t):=rtr^{-1} = t^{-1}$ or\/ $\sigma(t)=t g^k$ for some\/
   $k\in\ZZ$. In the latter case, one can change the
   generators of\/ $\cS(f)$ in such a way that the
   equation is satisfied with either\/ $k=0$ or\/ $k=1$.

   Finally, the following group structures for a reversible\/ $f$
   with involutory reversor\/ $r$ are possible after this reduction.
\begin{itemize}
% \item[\rm (1)] $\mathcal{R} (f) = \langle g\rangle \rtimes \langle
%    r\rangle \simeq C_{\infty}\rtimes C_{\infty}$, iff\/ $f$ has no
%    involutory reversors. In this case, one can choose\/ $t$ so that\/
%    $r^2=t$.
\item[\rm (1)] $\mathcal{R} (f) = \langle t\rangle \times
   \big(\langle g\rangle \rtimes \langle r\rangle\big)\simeq
   C_{\infty}\times D_{\infty}$, iff\/ $r$ commutes with\/ $t$. 
   In this case, also reversors of infinite order exist.
\item[\rm (2)] $\mathcal{R} (f) = \big(\langle t\rangle\times \langle
   g\rangle\big)\rtimes \langle r\rangle\simeq (C_{\infty}\times
   C_{\infty})\rtimes C_2$, iff either\/ $\sigma (t) = t^{-1}$
   $($which happens iff\/ {\em all} reversors are involutions\/$)$
   or\/ $\sigma(t)=tg$ $($in which case, once again, also reversors of
   infinite order exist\/$)$.
\end{itemize}
\end{theorem}
\begin{proof}
Since $rgr^{-1}$ is a non-trivial symmetry of $f$, we must have
$rgr^{-1}=t^{\varepsilon}g^{\ell}$ for some $\varepsilon,\ell\in \ZZ$,
not both $0$.  On the other hand,
$g^{-n}=f^{-1}=rfr^{-1}=rg^{n}r^{-1}=(rgr^{-1})^n =t^{n\varepsilon}
g^{n\ell}$, which implies $\varepsilon=0$ and $\ell=-1$. This shows
$rgr^{-1}=g^{-1}$. The statement about the order of $r$ is obvious
from the fact that $r^2$ is a symmetry.

Next, observe that $rtr^{-1}\neq 1$ is a symmetry, so that
$rtr^{-1}=t^{\varepsilon}g^k$ for some $\varepsilon,k\in\ZZ$, not both
$0$. Since $r^2$ commutes with $t$, one finds $t=r^2 t r^{-2} = r
t^{\varepsilon} g^k r^{-1} = (rtr^{-1})^{\varepsilon} rg^k r^{-1} =
(t^{\varepsilon} g^k)^{\varepsilon} g^{-k}=t^{\varepsilon^2}
g^{k(\varepsilon-1)}$.  This implies $\varepsilon^2=1$ and
$k(\varepsilon-1)=0$.  The solutions are $\varepsilon=-1$ together
with $k=0$, which means that $r$ is also a reversor for $t$, and
$\varepsilon=1$ together with an arbitrary $k\in \ZZ$, giving
$rtr^{-1}=t g^k$.

In the latter case, one may assume that $k\ge 0$ (otherwise,
replace the generator $g$ by $g^{-1}$). If $k>1$, one can
define a new generator $\tilde{t}=t g^{\lfloor k/2\rfloor}$, so that
$\tilde{t}$ and $g$ still generate the same group $\cS(f)$.
It is easy to check that this results in $r\tilde{t}r^{-1}
=\tilde{t}$ (resp.\ $\tilde{t} g$) depending on whether
$k$ was even (resp.\ odd).

For the final assertion, $\mathcal{R}(f) = \mathcal{S}(f) \rtimes
\langle r\rangle$ is clear by Lemma~\ref{semidirect}, where always
$\sigma(g)=g^{-1}$. The three cases now follow from the different
possibilities how $\sigma$ acts on $t$. If $r$ is a reversor for both
$t$ and $g$, we are again in the situation of part (1) of
Proposition~\ref{oneforall}. Otherwise, non-involutory reversors
exist, which must then be of infinite order.
\end{proof}

\begin{example} \cite{BRtorus}
Consider the matrices $M,R\in\mathrm{PGL}(4,\ZZ)$ given by
\[
    M \, = \, \left[ \begin{matrix}
    0 & 1 & 0 & 0 \\ 0 & 0 & 1 & 0 \\
    0 & 0 & 0 & 1 \\ -1 & 2 & 2 & 2 \end{matrix} \right]
\quad\mbox{and}\quad
    R \, = \, \left[ \begin{matrix}
    0 & 0 & 0 & 1 \\ 0 & 0 & 1 & 0 \\
    0 & 1 & 0 & 0 \\ 1 & 0 & 0 & 0 \end{matrix} \right] .
\]
It is easy to check that the involution $R$ conjugates $M$ into its
inverse.  As follows from \cite[Corollary~6]{BRtorus}, $M$ has
symmetry group $\mathcal{S}(M)\simeq C_{\infty}\times C_{\infty}$ in
$\mathrm{PGL}(4,\ZZ)$. One generator is $M$ itself, as this is a
matrix without roots in this matrix group, while the other generator
can be either chosen as
\[
    N \, = \, \left[ \begin{matrix}
    1 & 0 & -3 & 1 \\ -1 & 3 & 2 & -1 \\
    1 & -3 & 1 & 0 \\ 0 & 1 & -3 & 1 \end{matrix} \right]
\quad\mbox{or as}\quad
    N':= M\,N \, = \, \left[ \begin{matrix}
    -1 & 3 & 2 & -1 \\ 1 & -3 & 1 & 0 \\
    0 & 1 & -3 & 1 \\ -1 & 2 & 3 & -1 \end{matrix} \right] .
\]
\end{example}
Note that $N$ can neither possess a root in $\mathrm{GL} (4,\ZZ)$ nor
in $\mathrm{PGL} (4,\ZZ)$ because the sum of the square roots of the
eigenvalues of $N$ is not an integer.  The characteristic polynomial
of $N$ is $Q(x)=x^4-6x^3+22x^2-14x+1$, which is not self-reciprocal --
neither directly nor up to an overall sign. Consequently, $N$ is not
reversible within $G$, and neither within $\mathrm{GL} (4,\QQ)$,
compare \cite[Prop.~2]{BRtorus}, and the same statement applies to
$N'$. In fact, one quickly checks that $R$ and $N'$ commute, whence
$R':=R N'$ is a reversor of infinite order. The reversing symmetry
group thus has the structure $\mathcal{R} (M) = \langle N'\rangle
\times \big( \langle M\rangle\rtimes \langle R\rangle\big)\simeq
C_{\infty}\times D_{\infty}$, in line with case (1) of
Theorem~\ref{inf-inf}

\begin{remark}
The previous example can be considered within $\mathrm{GL}(4,\ZZ)$ as
well, i.e., as a toral automorphism. Note that the largest eigenvalue
of $M$ is a so-called Salem number. The characteristic polynomial of
$M$ is $P(x)=x^4-2x^3-2x^2-2x+1$, which is irreducible over $\ZZ$ (and
hence also over $\QQ$). Its roots are $\tau\pm\sqrt{\tau}$ (both real)
and $(1-\tau)\pm\sqrt{1-\tau}$ (both on the unit circle), where
$\tau=(\sqrt{5}+1)/2$ is the golden ratio. The symmetry is now
$\cS(M)\simeq C_2\times C_{\infty}\times C_{\infty}$, with
$C_2=\{\pm\one\}$, with the other details to be changed accordingly.
\end{remark}

\section{Comments and further directions}

Further extensions of the results along the lines of the previous
theorems are possible. In particular, one might want to extend the
setting to symmetry groups of the form $\cS(f)=C_m^{} \times
C_{\infty}^{\ell}$, with $\ell\ge 1$ and $m$ even, which occur for
toral automorphisms \cite{BRtorus} as a result of Dirichlet's unit
theorem. Since the methods should be clear from our above results, we
do not go into further detail.

Above, we have looked into the case that $\cS(f)=\cH\times\langle
g\rangle$ where $f$ was a power of $g$. In general, if $\cS(f)$ is a
finitely generated Abelian group, it is of the form $\cS(f)=\cF\times
C_{\infty}^k$ with $\cF$ a finite Abelian group, see
\cite[Thm.~I.8.5]{Lang}. For $f$ not of finite order, we might then
also assume that $f$ is an element of one of the $C_{\infty}$ factors.

However, in general, $\cS(f)$ need not be Abelian, whence there is no
compelling reason to start from a product structure such as $\cH\times
C_{\infty}$ (even with $\cH$ non-Abelian), as can be seen from the
possibility of $k$-symmetries. The general setting is then even more
involved, but can be handled by a computer assisted approach, e.g., as
in the classification of crystallographic point and space groups.

\smallskip
\section*{Acknowledgements}

M.B.\ would like to thank the School of Mathematics at UNSW
for hospitality, where most of this work was done. This work
was also supported by the German Research Council (DFG),
within the CRC 701, and by a 2005 UNSW Goldstar Award.

\end{document}